\documentclass[final,12pt]{article}
\usepackage{verbatim}
\usepackage{latexsym}
\usepackage{amsmath,amsthm,amsfonts}
\usepackage{eucal}
\usepackage[notref,notcite]{showkeys}

\usepackage{mathtext}
\usepackage[cp1251]{inputenc}
\usepackage[T2A]{fontenc}
\usepackage[dvips]{graphicx}
\usepackage{amsmath}
\usepackage{amssymb}
\usepackage{amsxtra}
\usepackage{latexsym}
\usepackage{ifthen}

\newtheorem{theorem}{Theorem}

\newtheorem{lemma}{Lemma}[section]

\newtheorem{corollary}{Corollary}[section]
\numberwithin{equation}{section}

\begin{document}

\title {\bf Existence and behavior of asymmetric traveling wave solutions
to thin film equation}

\author{{\normalsize\bf Roman M. Taranets}
\smallskip\\
{\normalsize\it Institute of Applied Mathematics and Mechanics of NAS of Ukraine},\\
{\normalsize\it R. Luxemburg str. 74, 83114 Donetsk, Ukraine}\\
{\normalsize\it taranets\_r@iamm.ac.donetsk.ua} }

\date{\today}
\maketitle

\setcounter{section}{0}

\begin{abstract}
We proved the existence and uniqueness of a traveling wave
solution to the thin film equation with a Navier slip condition at
the liquid-solid interface. We obtain explicit lower and upper
bounds for the solution and an absolute error estimate of
approximation of a solution to the thin films equation by the
traveling-wave solution.
\end{abstract}


\textbf{keywords:} {thin films, Navier-slip condition, asymmetric
traveling wave, lower and upper bounds for traveling waves,
absolute error estimate}

\section{Introduction}

The degenerate parabolic equation
\begin{equation}\label{0.0}
 h_t  + (h^n h_{xxx} )_x  = 0
\end{equation}
arises in description of the evolution of the height $y = h(t,x)$ of
a liquid film which spreads over a solid surface ($y = 0$) under the
action of the surface tension and viscosity in lubrication
approximation (see \cite{4,6}). Lubrication models have shown to be
extremely useful approximations to the full Navier-Stokes equations
for investigation of the thin liquid films dynamics, including the
motion and instabilities of their contact lines. For thicknesses in
the range of a few micrometers and larger, the choice of the
boundary condition at the solid substrate does not influence the
eventual appearance of instabilities, such as formation of fingers
at the three-phase contact line (see \cite{7,8}). For other
applications, such as for the dewetting of nano-scale thin polymer
film on a hydrophobic substrate the boundary condition at the
substrate appears to have crucial impact on the dynamics and
morphology of the film.

The exponent $n \in \mathbb{R}^{+}$ is related to the condition
imposed at the liquid-solid interface, for example, $n = 3$ for
no-slip condition, and $n \in (0,3)$ for slip condition in the
form
\begin{equation}\label{form}
v^{x} = \mu\,h^{n - 2} (v^{x})_{y} \text{ at } y = 0.
\end{equation}
Here, $v^x$ is the horizontal component of the velocity field,
$\mu$ is a non-negative slip parameter, and $\mu\,h^{n - 2}$ is
the weighted slip length. Distinguished are the cases $\mu = 0$
and $n = 2$. The first one corresponds to the assumption of a
no-slip condition, the second one to the assumption of a Navier
slip condition at the liquid-solid interface. The wetted region
$\{h > 0\}$ is unknown, hence the system is simulated as a
free-boundary problem, where the free boundary being given by
$\partial \{h > 0\}$, i.\,e. the triple junctions where liquid,
solid and air meet.

The main difficulty in studying equation (\ref{0.0}) is its
singular behaviour for $h = 0$. The mathematical study of
equation~(\ref{0.0}) was initiated by F.~Bernis, A.~Friedman
\cite{B8}. They showed the positivity property of solutions to
(\ref{0.0}) and proved the existence of nonnegative generalized
solutions of initial--boundary problem with an arbitrary
nonnegative initial function from $H^1$. More regular
(\emph{strong} or \emph{entropy}) solutions have been constructed
in \cite{B2,B14}. One outstanding question is whether zeros
develop in finite time, starting with a regular initial data. What
is known is that with periodic boundary conditions, for $n
\geqslant 3.5$ this does not occur \cite{B2,B8}, while for $n <
3/2$ the solution develops zeroes in a finite time \cite{10}. One
way of looking at the problem (\ref{0.1}) has been to study
similarity solutions to (\ref{0.0}) in the form $h(x,t) =
t^{-\alpha} H(x\,t^{- \beta})$, where $n\alpha + 4\beta =1$ (see
\cite{9}). In the paper \cite{1}, the authors proved also
existence dipole solutions and found their asymptotic behaviour.
We note that the solutions of such type do not exist in the case
$n \geqslant 2$, however, there exists a traveling wave solution
(see \cite{4}).

In the present paper, we concentrate on a traveling wave solution
to (\ref{0.0}) at $n = 2$, namely, we consider the following
problem with a regular initial function:
\begin{equation}\label{0.1}
\left\{ \begin{array}{l}
h_t  + (h^2 h_{xxx} )_x  = 0, \\
h(s_1 ) = h(s_2 ) = 0, \\
h_x (s_1 ) = \theta  > 0,\  h^2 h_{xxx}  = 0 \text{ at } x = s_1.\\
 \end{array} \right.
\end{equation}
System (\ref{0.1}) describes the growth of dewetted regions in the
film. Fluid transported out of the growing dry regions collects in a
ridge profile which advances into the undisturbed fluid (see
Figure~1). Under ideal conditions, it could be imagined that dry
spots could grow indefinitely large. By conservation of mass, the
growing holes would shift fluid into the ever-growing rims. In our
situation, large length scale to limit the sizes of these structures
is absent, and we might expect the motion and growth of the ridges
to approach scale-invariant self-similar form. At the same time, the
ridge profiles have a pronounced asymmetry (see \cite{3}).

\begin{figure}[t]
\includegraphics[width= \textwidth, keepaspectratio]{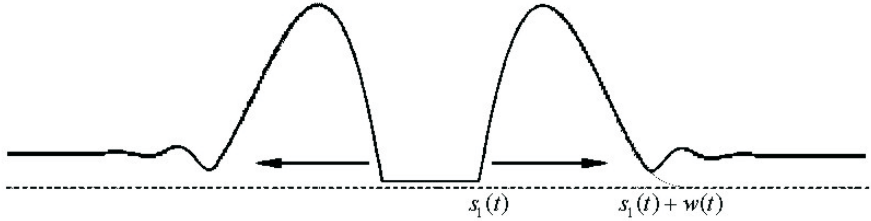}
\caption[1]{Sketch of the cross-section of a dewetting film after
rupture, showing the expanding dewetted residual layer (also
called a ''hole'' or ''dry spot'') in the middle and the adjacent
''dewetting ridges'' moving into the surrounding undisturbed
uniform film \cite{3}.}
\end{figure}


In the problem at hand, $x = s_1(t)$ is the position of the former
moving interface, i.\,e. the contact line, while the position of
the latter interface will give an effective measure of the width
of the ridge, $x = s_2(t) = s_1(t) + w(t)$. The ridge is assumed
to be moving forward, $\dot{s}_1(t) > 0$, corresponding to an
expanding hole. The arbitrary positive parameter $\theta$
corresponds to the contact angle of liquid-solid interface. Thus,
we can control a dewetted region in the film by the contact line
and obtain asymmetry profiles of solutions. As the paper \cite{3}
has shown that the axisymmetric profile can be analyzed within a
one-dimensional thin-film model. The authors found matched
asymptotic expansion, speed and structure of the profile, in
particular, they obtained that
\begin{equation}\label{as}
h(x,t) \sim A\, \dot{s_1}^{1/2} (s_2  - x)^{3/2} \text{ at } x =
s_2,
\end{equation}
where the asymptotic constant $A = 2(2/3)^{1/2}$.

Hereinafter, we assume that the contact line moves with a constant
velocity (${\rm{v}}$) and the width of the ridge ($w$) is a
constant. As in \cite{3}, we are going to look for a solution to
(\ref{0.1}) in the following form
$$
h(x,t) = \hat{h}(\xi ),  \text{ where } \xi  = x - {\rm{v}}t,\
s_1 = {\rm{v}}t,\  s_2  = {\rm{v}}t + w.
$$
We can remove v and $\theta$ from the resulted problem by
rescaling appropriately,
$$
\hat{h} = \tfrac{\theta ^3}{\rm{v}}\varphi ,\ \xi  = \tfrac{\theta
^2 }{{\rm{v}}}\eta ,\  w = \tfrac{\theta ^2 }{{\rm{v}}}d.
$$
As a result, we obtain the following problem for the traveling-wave
\begin{equation}\label{0.2}
\left\{ \begin{array}{l}
\varphi (\eta) \varphi'''(\eta)  = 1,\  \varphi (\eta ) \geqslant 0, \\
\varphi (0) = \varphi (d) = 0,\ \varphi'  (0) = 1,\  0 \leqslant \eta  \leqslant d. \\
\end{array} \right.
\end{equation}
Boatto et al. \cite{4} reduced this problem to the problem of
finding a co-dimension one orbit of a second-order ODE system
connecting equilibria. Hence generically solutions will exist but
only for isolated values of the free parameter $d$. The parameter
$d$ was found in \cite{3} by integration (\ref{0.2}), and
\begin{equation}\label{d}
    d=1/2.
\end{equation}

Our paper is organized as follows. In Section~2 we prove the
existence and uniqueness of a traveling wave solution to the
problem (\ref{0.2}) (Theorem~\ref{Th1}). Lower and upper bounds
for the traveling wave solution are contained in Section~3
(Theorem~\ref{Th2}). We note that the bounds assert that the
constant $A$ of (\ref{as}) must be from the interval
$[4\sqrt{2}/3, 4 \sqrt{6}/3]$ (Corollary~3.1). In Section~4 we
find an absolute error estimate of approximation of a solution to
(\ref{0.1}) by the traveling-wave solution (Theorem~\ref{Th3}).

\section{Existence of the traveling wave solution}

Below we prove the existence and uniqueness traveling wave solution
to the problem (\ref{0.2}). Our proof is based on some modification
of the proof of the existence and uniqueness dipole solutions from
\cite{1}.

\begin{theorem}\label{Th1}
There exists a unique solution $\varphi (\eta )$ to the problem
(\ref{0.2}) such that $\varphi (\eta ) \in C^3 (0,d) \cap C^1
[0,d]$ and $\varphi (\eta ) > 0$ for $0 < \eta  < d$.
\end{theorem}

\noindent  First, we prove the following auxiliary lemma:

\begin{lemma}\label{l1.2}
Assume that $\varphi  \in C^3 (0,d) \cap C^1 [0,d]$, $ \varphi (0) =
\varphi (d) = 0,\ \varphi'  (0) = 1$, $\varphi
> 0$ and $\varphi''' (\eta )\geqslant 0$ in $(0,d)$. Then
$\varphi$ has a unique maximum and $\varphi'  (d) \leqslant 0$.
\end{lemma}

\begin{proof}
Since $\varphi''' (\eta ) \geqslant 0$ we have that $ \varphi'(\eta
)$ is convex, $\varphi''(\eta )$ is increasing. By Rolle's theorem
$\varphi '(\eta )$ has at least one zero in $(0,d)$, and
$\varphi'(\eta )$ has no more than two zeroes in $(0,d)$ by
convexity. Let $d_1 ,d_2  \in (0,d) : \varphi'  (d_1 ) = \varphi'
(d_2 ) = 0 $. Then, by Rolle's theorem, there exists $d_3 \in (d_1
,d_2 ): \varphi'' (d_3 ) = 0$ whence $\varphi (d) \geqslant 0$. In
view of $\varphi  > 0$ in $ (0,d)$ and $\varphi (d) = 0$, we obtain
$d_2 = d$ and $\varphi' (d) = 0$. This proves that $\varphi'(\eta )$
has exactly one zero in $(0,d)$ and hence $\varphi (\eta )$ has a
unique maximum. Now $ \varphi' (d) < 0$ follows easily.
\end{proof}

\begin{proof}[Proof of Theorem~\ref{Th1}] \emph{Green's function.} We define a Green's
function $G(\eta ,t)$ by
\begin{equation}\label{1.1}
\left\{ \begin{array}{l}
G'''(\eta ,t) = \delta (\eta  - t),\  0 \leqslant \eta  \leqslant d,\ 0 \leqslant t \leqslant d, \\
G(0,t) = G(d,0) = G'(0,t) = 0,\  0 \leqslant t \leqslant d,\\
\end{array} \right.
\end{equation}
where $d = 1/2$. By explicit computation, we find that
\begin{equation}\label{1.2}
0 \leqslant G(\eta ,t) = \left\{ \begin{array}{l}
 2(t - d)^2 \eta ^2 \text{ if } 0 \leqslant \eta  \leqslant t \leqslant d, \\
 2(t - d)^2 \eta ^2  - d(\eta  - t)^2 \text{ if } 0 \leqslant t \leqslant \eta  \leqslant d, \\
 \end{array} \right.
\end{equation}
whence
\begin{equation}\label{1.2*}
\int\limits_0^\eta  {G(\eta ,t)\,dt}  = \tfrac{2}{3}\eta ^3 (d -
\eta )(1 - \eta ),\ \int\limits_\eta ^d {G(\eta ,t)\,dt}  =
\tfrac{2}{3}\eta ^2 (d - \eta )^3
\end{equation}
if $0 < \eta  < d$, and
\begin{equation}\label{1.3}
G(\eta ,t) \leqslant C\,t^2 (d - t),\quad |G'(\eta ,t)| \leqslant
C\,t(d - t)\ \forall\, \eta ,\,t \in [0,d].
\end{equation}

\medskip

\noindent \emph{Approximating problems.} For each positive integer
$k$ we consider the problem
\begin{equation}\label{1.4}
\left\{ \begin{array}{l}
 \varphi_{k}''' (\eta ) = \varphi _k^{ - 1} \text{ for }  0 < \eta  < d, \\
 \varphi _k (0) = \varphi _k (d) = \frac{1}{k},\  \varphi' _{k} (0) = 1, \\
 \varphi _k (\eta ) \in C^3 [0,d],\  \varphi _k (\eta ) > 0   \text{ for } 0
 \leqslant \eta  \leqslant d. \\
 \end{array} \right.
\end{equation}
Consider the closed convex set
$$
S = \{ v \in C[0,d]:v \geqslant 1/k \text{ in } (0,d)\}
$$
and the nonlinear operator $\Phi _k$ defined by
$$
\Phi _k v(\eta ) = \tfrac{1}{k} + 2\eta (d - \eta ) +
\int\limits_0^d {G(\eta ,t)v^{ - 1} (t)\,dt},
$$
where $G(\eta ,t)$ is from (\ref{1.2}). The operator $\Phi _k$
mapping $S$ into $S$ is continuous. Moreover, $\Phi _k (S)$ is (for
each $k$) a bounded subset of $C^3 [0,d]$ and hence a relatively
compact subset of $S$. By Schauder's fixed-point theorem, there
exists $\varphi _k \in \Phi _k (S)$ such that $\Phi _k \varphi _k  =
\varphi _k$. This is the desired solution of the problem~(\ref{1.4}).
Note that $\varphi _k$ satisfies
\begin{equation}\label{1.5}
\begin{gathered}
\varphi _k (\eta ) = \tfrac{1}{k} + 2\eta (d - \eta ) +
\int\limits_0^d {G(\eta ,t)\varphi''' _{k} (t)\,dt} ,\hfill\\
\varphi' _{k} (\eta ) = 2(d - 2\eta ) + \int\limits_0^d {G'(\eta
,t)\varphi''' _{k} (t)\,dt}. \hfill
\end{gathered}
\end{equation}
In view of Lemma~\ref{l1.2} (applied to $\varphi _k  - 1/k$),
there exists a unique point $m_k$ in which the maximum of $\varphi
_k$ is attained. Therefore,
$$
\varphi _k (\eta ) \nearrow  \text{ in } (0,m_k )\ \text{and}\
 \varphi _k (\eta ) \searrow \text{ in } (m_k ,d).
$$

\medskip

\noindent \emph{Estimates.} Since $G \geqslant 0$ we get
$$
\begin{array}{r}
\varphi _k (\eta ) - \frac{1}{k} = 2\eta (d - \eta ) +
\int\limits_0^d {G(\eta ,t)\varphi _k^{ - 1} (t)\,dt}  \geqslant
 \varphi _k^{ - 1} (\eta )\int\limits_0^\eta  {G(\eta ,t)\,dt}  \\
= \tfrac{2}{3} \eta ^3 (d - \eta )(1 - \eta ) \varphi _k^{ - 1} (\eta )
\geqslant \frac{2}{3}  \eta ^3 (d - \eta )^2\varphi _k^{ - 1} (\eta ),  \\
\end{array}
$$
whence
$$
\varphi _k (\eta ) \geqslant \sqrt {\tfrac{2}{3}} \eta ^{3/2} (d -
\eta ) \geqslant C\,\eta (d - \eta )^{3/2} \quad {\rm{if}}\quad 0
< \eta  < m_k,
$$
and
$$
\begin{array}{r}
\varphi _k (\eta ) - \frac{1}{k} = 2\eta (d - \eta ) +
\int\limits_0^d {G(\eta ,t)\varphi _k^{ - 1} (t)\,dt}  \geqslant
\varphi _k^{ - 1} (\eta )\int\limits_\eta ^d {G(\eta ,t)\,dt}  \\
= \tfrac{2}{3} \eta ^2 (d - \eta )^3 \varphi _k^{ - 1} (\eta ), \\
\end{array}
$$
whence
$$
\varphi _k (\eta ) \geqslant \sqrt {\tfrac{2}{3}} \eta (d - \eta
)^{3/2} \quad {\rm{if}}\quad m_k  < \eta  < d.
$$
Hence,
\begin{equation}\label{1.6}
\varphi _k (\eta ) \geqslant C\,\eta (d - \eta )^{3/2} \quad
\forall\, \eta  \in (0,d).
\end{equation}
Next we deduce from the differential equation that, for all $\eta
\in (0,d)$,
\begin{equation}\label{1.7}
\eta (d - \eta )\varphi''' _{k} (\eta ) = \eta (d - \eta )\varphi
_k^{ - 1} (\eta )  \mathop \leqslant \limits^{(\ref{1.6})} C\,(d -
\eta )^{ - 1/2} ,
\end{equation}
where the right-hand side is an integrable function.

\medskip

\noindent \emph{Passing to the limit.} From (\ref{1.5}),
(\ref{1.3}) and (\ref{1.7}) it follows that $\varphi _k$ is
bounded in $C^1 [0,d]$. Therefore, there exists a subsequence,
again denoted by $\varphi _k$, and a function $\varphi$ such that
$\varphi _k \mathop  \to \varphi$ uniformly on $[0,d]$ as $ k \to
\infty $. Thus $\varphi (0) = \varphi (d) = 0$ and, by
(\ref{1.6}), $\varphi
> 0$ in $(0,d)$. Hence, for each compact subset $I$ of $(0,d)$ we
have
$$
\varphi''' _{k} (\eta ) = \varphi _k^{ - 1} (\eta )\mathop  \to
\limits_{k \to \infty } \varphi ^{ - 1} (\eta ) \text{ in } C(I),
$$
and, by (\ref{1.7}) and Lebesgue's dominated convergence theorem,
\begin{equation}\label{1.8}
\eta (d - \eta )\varphi''' _{k} (\eta )\mathop \to \limits_{k \to
\infty } \eta (d - \eta )\varphi ^{ - 1} (\eta )\text{ in } L^1
(0,d).
\end{equation}
Since $\varphi''' _{k} \mathop  \to \limits_{k \to \infty }
\varphi''' $ in the distribution sense, it follows that
\begin{equation}\label{1.9}
\varphi ''' (\eta ) = \varphi ^{ - 1} (\eta )\text{ in } (0,d),
\end{equation}
i.e. $\varphi$ satisfies the differential equation. Moreover, from
(\ref{1.8}), (\ref{1.9}), (\ref{1.5}) and (\ref{1.3}) we deduce that
$\varphi _k \mathop  \to  \varphi$ in $C^1 [0,d]$ as $k \to \infty
$, and hence $\varphi$ also satisfies $\varphi' (0) = 1$. This
completes the proof of the existence.

\medskip

\noindent \emph{Uniqueness}. Let $\varphi_1$ and $\varphi_2$ be
two solutions of the problem (\ref{0.2}) and set $v = \varphi_1 -
\varphi_2$. Since  $v''' = \varphi_1^{-1} - \varphi_2^{-1}$ and
the function $v \mapsto v^{-1}$ is decreasing we deduce that $v\,v
''' \leqslant 0$. Since $v\,v''' \leqslant 0$ and $v(0) = v(d) = v
'(0) = 0$, we conclude that $v \equiv 0$. This completes the proof of
Theorem~\ref{Th1}.
\end{proof}

\section{Lower and upper bounds for the traveling wave solution}

Integrating (\ref{0.2}) with respect to $\eta$, we arrive at the
following problem
\begin{equation}\label{2.1}
\left\{ \begin{array}{l}
 \varphi(\eta) \varphi''(\eta)  = \tfrac{1}{2}(\varphi'(\eta) )^2
 + \eta  - d,\  \varphi (\eta ) \geqslant 0, \\
 \varphi (0) = 0,\,\varphi'  (0) = 1,\  0 \leqslant \eta  \leqslant d, \\
 \end{array} \right.
\end{equation}
where $d$ is from (\ref{d}). Analyzing the behaviour of a solution
to (\ref{2.1}), we find explicit lower and upper bounds for the
solution.

\begin{theorem}\label{Th2}
Let $\varphi (\eta )$ be a solution from Theorem~\ref{Th1}. Then the
following estimates are valid
\begin{equation}\label{2.5}
\varphi _{\min} (\eta ) \leqslant \varphi (\eta ) \leqslant
\varphi _{\max} (\eta ) \Leftrightarrow \tfrac{4\sqrt{2}}{3}\eta
(d - \eta )^{3/2}\leqslant \varphi (\eta ) \leqslant
\tfrac{{4\sqrt 6 }}{3}\,\eta (d - \eta )^{3/2}
\end{equation}
for all $\eta  \in [0,d]$ $($see Figure~2$)$.
\end{theorem}

\begin{corollary}
In particular, from (\ref{2.5}) it follows that
\begin{equation}\label{2.5'}
\tfrac{4\sqrt{2}}{3} \dot{s}_1^{1/2} \tfrac{x - s_1}{s_2 - s_1}
(s_2 - x)^{3/2} \leqslant h(x,t) \leqslant \tfrac{4\sqrt 6}{3}
\dot{s}_1^{1/2} \tfrac{x - s_1}{s_2 - s_1} (s_2 - x)^{3/2}
\end{equation}
for all $x \in [s_1, s_2]$.
\end{corollary}

\begin{figure}[t]
\includegraphics[width= \textwidth, keepaspectratio]{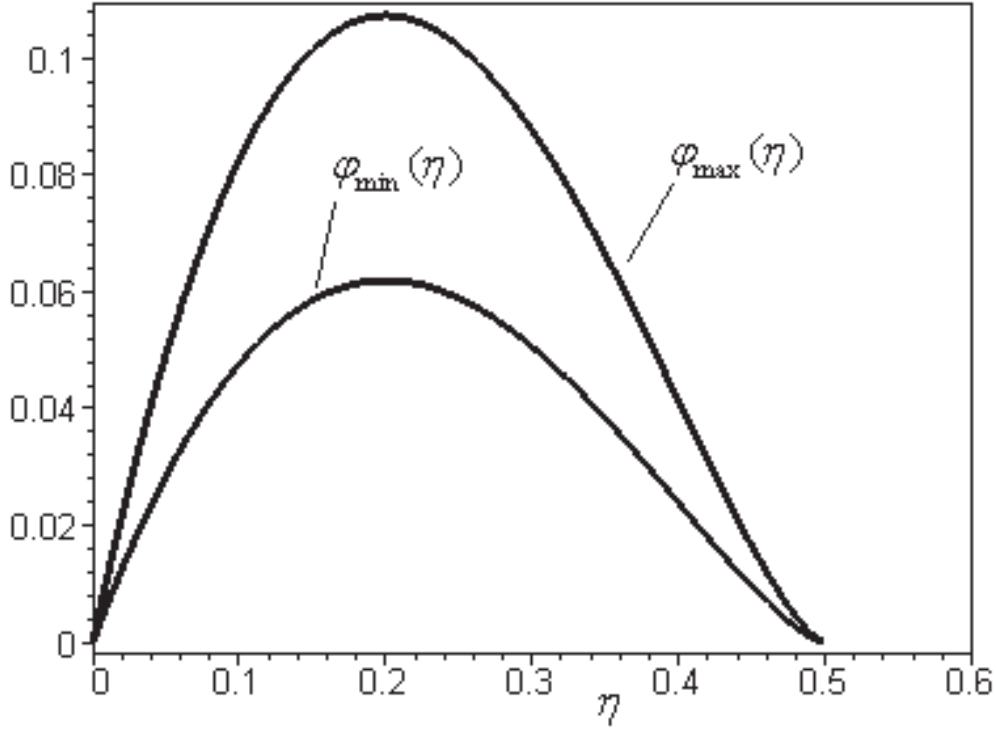}
\caption[2]{Lower and upper bounds for the traveling wave
solution}
\end{figure}

\begin{lemma}\label{l2.1}
The function $ \varphi_0 (\eta ) = A_0\,\eta (d - \eta )^{3/2}$
$(A_0  > 0)$ satisfies the inequalities
\begin{equation}\label{2.2}
\varphi \varphi''  \geqslant \tfrac{1}{2}(\varphi' )^2  + \eta - d
\  \forall\, \eta  \in [0,d] \text{ if } A_0^2 d^2 \leqslant 5/3,
\end{equation}
\begin{equation}\label{2.3}
\varphi \varphi''  \leqslant \tfrac{1}{2}(\varphi')^2  + \eta  - d
\  \forall\, \eta  \in [0,d]  \text{ if } A_0^2 d^2 \geqslant 8/3.
\end{equation}
\end{lemma}

\begin{proof}
Indeed, the function $\varphi_0 (\eta )$ satisfies the equation
$$
\varphi _0 \varphi'' _{0}  = \tfrac{1}{2}(\varphi'_{0} )^2 + f
(\eta ),
$$
where
$$
f(\eta ): = \tfrac{5}{8}A_0^2 (d - \eta )\left( {\eta  -
\tfrac{2d(1 - \sqrt 6 )}{5}} \right)\left( {\eta  - \tfrac{2d(1 +
\sqrt 6 )}{5}} \right) \leqslant 0 \ \forall\, \eta  \in [0,d].
$$
From
$$
f(\eta ) \geqslant \eta  - d \Leftrightarrow 5\eta ^2  - 4d\eta  -
4d^2  + \tfrac{8}{A_0^2} \geqslant 0 \  \forall \,\eta  \in [0,d],
$$
$$
D = 4d^2  + 20d^2  - \tfrac{40}{A_0^2} = \tfrac{24}{A_0^2} \left(
{A_0^2 d^2  - \tfrac{5}{3}} \right) \leqslant 0 \text{ if } A_0^2
d^2  \leqslant 5/3
$$
we obtain (\ref{2.2}). In a similar way, we obtain (\ref{2.3}) for
$f(\eta ) \leqslant \eta  - d$ $\forall\, \eta  \in [0,d]$ if
$A_0^2 d^2  \geqslant 8/3$.
\end{proof}

\begin{lemma}\label{l2.2}
The function $\varphi_{\min} (\eta ) = A_1 \eta (d - \eta )^{3/2}$
$(A_1  > 0)$ if $A_1^2 d^2  \leqslant 8/9$ is a lower bound for
the solution $\varphi (\eta )$ of (\ref{2.1})$:$
$$
\varphi _{\min}(\eta ) \leqslant \varphi (\eta ) \ (\text{i.\,e. }
\varphi _{\min} (\eta ) - \varphi (\eta ) \leqslant 0) \ \forall\,
\eta  \in [0,d].
$$
\end{lemma}

\begin{proof} ''Contraction Principle''. Let us define
$v(\eta):= \varphi _{\min}(\eta) - \varphi (\eta)$. Suppose that
there exists a point $\eta _0  \in [0,d]$ such that $ v(\eta_0) >
0$ then $\eta _0$ is a point of maximum for $v(\eta)$, i.\,e.
$v'(\eta _0 ) = 0 \Leftrightarrow \varphi '_{\min} (\eta _0 ) =
\varphi '(\eta _0 ) = M$ and $v''(\eta _0) < 0$. From (\ref{2.1})
and (\ref{2.2}) we deduce that
$$
\left\{ \begin{array}{l}
\varphi \varphi''  = \tfrac{1}{2}(\varphi')^2  + \eta  - d \\
\varphi _{\min} \varphi'' _{\min} \geqslant
\tfrac{1}{2}(\varphi' _{\min})^2  + \eta  - d \\
\end{array} \right.\Rightarrow
$$
$$
v''(\eta) \geqslant \tfrac{1}{2}\left( {\tfrac{(\varphi
_{\min,\eta } )^2 }{\varphi _{\min}} - \tfrac{(\varphi _\eta  )^2
}{\varphi }} \right) + (\eta - d)\left({\tfrac{1}{\varphi _{\min}}
- \tfrac{1}{\varphi }} \right),
$$
whence
$$
\underbrace {v''(\eta_0)}_{ < 0} \geqslant \underbrace
{\tfrac{-v(\eta _0)}{\varphi _{\min} (\eta _0 )\varphi (\eta _0
)}}_{ < 0}\underbrace {\left( {\tfrac{1}{2}M^2  + \eta _0 - d}
\right)}_{?}.
$$
Using $M = A_1 (d - \eta _0 )^{1/2} (d - 5\eta _0 /2)$, we find
$$
\tfrac{1}{2}M^2  + \eta _0  - d = (d - \eta _0 )\bigl[
\tfrac{1}{2} A_1^2 (d - 5\eta _0 /2)^2  - 1 \bigr] \leqslant 0 \
\forall \,\eta _0  \in [0,d] \text{ if } A_1^2 d^2 \leqslant 8/9.
$$
Thus we obtain a contradiction with our assumption, which proves
the assertion of
Lemma~\ref{l2.2}.
\end{proof}

\begin{lemma}\label{l2.3}
The function $ \varphi _{\max} (\eta ) = A_2\, \eta (d - \eta
)^{3/2}$ $(A_2  > 0)$ if $A_2^2 d^2  \geqslant 8/3$ is an upper
bound for the solution $\varphi (\eta )$ of (\ref{2.1})$:$
$$
\varphi (\eta ) \leqslant \varphi _{\max} (\eta ) \ (\text{i.\,e.
} \varphi (\eta ) - \varphi _{\max} (\eta ) \leqslant 0) \ \forall
\eta  \in [0,d].
$$
\end{lemma}

\begin{proof}
''Contraction Principle''. Let us define $v(\eta):= \varphi (\eta)
- \varphi _{\max}(\eta)$. Suppose that there exists a point $\eta
_0  \in [0,d]$ such that $ v(\eta_0)
> 0$ then $\eta _0$ is a point of maximum for $v(\eta)$, i.\,e.
$v'(\eta _0 ) = 0 \Leftrightarrow \varphi '_{\max} (\eta _0 ) =
\varphi '(\eta _0 )$ and $v''(\eta _0) < 0$. Moreover, from
$\varphi ''_{\max} (\eta ) = 3A_1 (d - \eta )^{ - 1/2} [5\eta /4 -
d]$ we find that $\varphi ''_{\max} (\eta ) \leqslant 0$ for all
$\eta \in [0,0.4]$ and $\varphi ''_{\max} (\eta ) > 0$ for all
$\eta \in (0.4,d]$. From (\ref{2.1}) and (\ref{2.3}) we deduce that
\begin{equation}\label{2.4}
\left\{ \begin{array}{l}
\varphi \varphi''  = \tfrac{1}{2}(\varphi' )^2  + \eta  - d \\
\varphi _{\max} \varphi'' _{\max}
\leqslant \tfrac{1}{2}(\varphi' _{\max})^2  + \eta  - d\\
\end{array} \right.
\Rightarrow \varphi(\eta _0 ) \varphi'' (\eta _0 )  - \varphi
_{\max}(\eta _0 ) \varphi'' _{\max}(\eta _0 ) \geqslant 0,
\end{equation}
whence
$$
\underbrace{\varphi (\eta _0 )\varphi''(\eta _0 ) - \varphi
_{\max} (\eta _0 )\,\varphi'' _{\max} (\eta _0 )}_{\geqslant 0} =
\underbrace {\varphi (\eta _0 ) v''(\eta _0 )}_{ \leqslant 0} +
\underbrace {\varphi''_{\max} (\eta _0 )\, v(\eta _0 )}_{
\leqslant 0}.
$$
This contradicts to our assumption if $\eta _0  \in
[0,0.4]$.

Now, let $\eta _0  \in (0.4,d]$. In this case, if $\varphi ''(\eta
_0 ) \leqslant 0$, and we obtain a contradiction immediately
from (\ref{2.4}). If $\varphi ''(\eta _0 )
> 0$ then we rewrite (\ref{2.4}) in equivalent form:
\begin{equation}\label{2.04}
\left\{ \begin{array}{l}
\tfrac{1}{2}(\varphi^2 )''  = \frac{3}{2}(\varphi')^2  + \eta  - d \\
\tfrac{1}{2}(\varphi_{\max}^2)''
\leqslant \tfrac{3}{2}(\varphi' _{\max} )^2  + \eta  - d \\
\end{array} \right.
\Rightarrow (\varphi^2(\eta _0 )  - \varphi_{\max}^2(\eta _0 ))''
\geqslant 0,
\end{equation}
whence
$$
\underbrace{(\varphi^2(\eta _0 )  - \varphi_{\max}^2(\eta _0
))''}_{\geqslant 0} = \underbrace{(\varphi''(\eta _0 )  +
\varphi''_{\max}(\eta _0))}_{>0} \underbrace{v''(\eta _0 )}_{<0}
$$
and we arrive at a contradiction. Thus, Lemma~\ref{l2.3} is
proved.
\end{proof}

As a result of Lemmata~\ref{l2.2} and \ref{l2.3}, we obtain lower
(more exact in comparison with (\ref{1.6})) and upper bounds
for the solution of (\ref{2.1}), and consequently for the solution
of (\ref{0.2}).

\section{Absolute error estimate of approximation of a solution by a traveling-wave solution}

The next theorem contains an absolute error estimate of
approximation of a solution (e.g., a generalized solution) by a
traveling-wave solution.

\begin{theorem}\label{Th3}
Let $h(x,t)$ be a solution and $\hat{h}(\xi )$ be the
traveling-wave solution to the problem (\ref{0.1}). Then the
following estimates hold
\begin{equation}\label{3.1}
\begin{array}{l}
\mathop {\sup }\limits_{x \in [s_1 ,s_2 ]} | {h(x,t) - \hat{h}(\xi
)} | \leqslant \tfrac{\sqrt{6}}{6}\theta (s_2 - s_1)^{1/2}  \text{ if } |h_x (s_2 )| > \theta , \\
\mathop {\sup }\limits_{x \in [s_1 ,s_2 ]} \!\!\! \bigl|h(x,t) -
\hat{h}(\xi )  - 2 \theta s_1^{1/2} \!(s_2 - s_1)^{1/2} \bigr|
\!\leqslant \!\tfrac{\sqrt{6}}{6}\theta (s_2 - s_1)^{1/2}
\text{ if } |h_x (s_2 )| \leqslant \theta, \\
\end{array}
\end{equation}
where $\xi  = x - {\rm{v}}t,\ s_1  = {\rm{v}}t$ and $s_2 = s_1  +
w$.
\end{theorem}

\begin{proof}[Proof of Theorem~\ref{Th3}] We make  the following
change of variables in (\ref{0.1})
$$
h(x,t) \mapsto f(\xi ,t), \text{ where }  \xi  = x - {\rm{v}}t.
$$
As a result, we obtain the following problem
\begin{equation}\label{3.2}
\left\{ \begin{array}{l}
 f_t  - {\rm{v}}\,f_\xi   + (f^2 f_{\xi \xi \xi } )_\xi   = 0, \\
 f(0,t) = f(w,t) = 0,\quad f^2 f_{\xi \xi \xi }  = 0 \text{ at }\xi  = 0, \\
 f_\xi  (0,t) = \theta  > 0. \\
 \end{array} \right.
\end{equation}
Multiplying (\ref{3.2}$_1$) by $- f_{\xi \xi } (\xi ,t)$ and
integrating with respect to $\xi$, we get
$$
\begin{array}{c}
\tfrac{1}{2} \tfrac{d}{dt} \int\limits_0^w {f_\xi ^2 (\xi
,t)\,d\xi } =  - {\rm{v}}\int\limits_0^w {f_\xi f_{\xi \xi }\,
d\xi } + \int\limits_0^w {(f^2 f_{\xi \xi \xi } )_\xi f_{\xi \xi
}\, d\xi } = - \tfrac{ {\rm{v}}}{2} \int\limits_0^w
{ \tfrac{\partial}{\partial \xi}(f_\xi ^2 )\,d\xi }  +  \\
+ \int\limits_0^w {(f^2 f_{\xi \xi \xi } )_\xi  f_{\xi \xi }
\,d\xi }  \mathop = \limits^{(\ref{3.2}_2),(\ref{3.2}_3)}
\tfrac{\rm{v}}{2} (\theta ^2  - f_\xi ^2 (w,t)) - \int\limits_0^w
{f^2 f_{\xi \xi \xi }^2\, d\xi },
\end{array}
$$
whence
\begin{multline}\label{3.3}
\tfrac{1}{2} \tfrac{d}{dt}\int\limits_0^w {f_\xi ^2 (\xi ,t)\,d\xi
} + \int\limits_0^w {f^2 f_{\xi \xi \xi }^2 \,d\xi }  =
\tfrac{\rm{v}}{2}(\theta ^2  - f_\xi ^2 (w,t)) \Rightarrow \\
\tfrac{d}{dt} \int\limits_0^w {f_\xi ^2 (\xi ,t)\,d\xi } \leqslant
{\rm{v}}(\theta ^2  - f_\xi ^2 (w,t)).
\end{multline}
Integrating (\ref{3.3}) with respect to time, we find
\begin{equation}\label{3.4}
\| {f_\xi  (\xi ,t)} \|_{L^2 (0,w)}^2  \leqslant \| {f_\xi (\xi
,0)} \|_{L^2 (0,w)}^2  + {\rm{v}}\int\limits_0^t {(\theta ^2 -
f_\xi ^2 (w,t))\,dt} .
\end{equation}
From (\ref{3.4}) it follows that
\begin{equation}\label{3.4*}
\left\| {f_\xi  (\xi ,t)} \right\|_{L^2 (0,w)}^2  \leqslant
\left\{
\begin{array}{c}
\left\| {f_\xi  (\xi ,0)} \right\|_{L^2 (0,w)}^2 \text{ if } |f_\xi  (w,t)| > \theta , \\
\left\| {f_\xi  (\xi ,0)} \right\|_{L^2 (0,w)}^2  + {\rm{
v}}\theta ^2 t \text{ if }|f_\xi  (w,t)| \leqslant \theta . \\
\end{array} \right.
\end{equation}
From this, by virtue of uniqueness of the traveling-wave solution
$\hat{h}(\xi)$ (see Theorem~\ref{Th1}), $f(\xi ,0) = \hat{h}(\xi)
$ we deduce from (\ref{3.4*}) that
$$
\begin{array}{l}
\| {f_\xi  (\xi ,t) - \hat{h}_\xi  (\xi )} \|_{L^2 (0,w)}
\leqslant
2\| {\hat{h}_\xi  (\xi )} \|_{L^2 (0,w)} \text{ if } |f_\xi  (w,t)| > \theta , \\
\| {f_\xi  (\xi ,t) - \hat{h}_\xi  (\xi )} \|_{L^2 (0,w)}
\leqslant 2\| {\hat{h}_\xi  (\xi )} \|_{L^2 (0,w)}  + 2\theta
\sqrt
{{\rm{v}}\,t} \text{ if } |f_\xi  (w,t)| \leqslant \theta . \\
\end{array}
$$
Therefore taking into account the embedding $\mathop {H^1}
\limits^{o}(0,w) \subset C[0,w]$, we find that
\begin{equation}\label{3.06}
\begin{gathered}
\mathop {\sup }\limits_{\xi  \in [0,w]} | {f(\xi ,t) - \hat{h}(\xi
)} | \leqslant 2 w^{1/2}\| {\hat{h}_\xi  (\xi )} \|_{L^2 (0,w)} \text{ if } |f_\xi  (w,t)| > \theta , \hfill\\
\mathop {\sup }\limits_{\xi  \in [0,w]} |f(\xi ,t) - \hat{h}(\xi )
| \leqslant 2 w^{1/2} (\| {\hat{h}_\xi  (\xi )} \|_{L^2 (0,w)}+
\theta \sqrt {{\rm{v}}\,t}) \text{ if } |f_\xi  (w,t)| \leqslant
\theta .\hfill
\end{gathered}
\end{equation}
Since $\hat{h}_\xi  (\xi ) = \theta \varphi'_{\eta}(\eta)$
($\varphi(\eta)$ is from Theorem~\ref{Th1}) and $\varphi(\eta)$
has a unique maximum in $(0,d)$, due to (\ref{2.5}) we conclude
that
\begin{equation}\label{3.7}
\| {\hat{h}_\xi  (\xi )} \|_{L^2 (0,w)} \leqslant \tfrac{
\sqrt{6}}{12}\theta.
\end{equation}
Thus, we  from (\ref{3.06}) and (\ref{3.7}) arrive at
\begin{equation}\label{3.6}
\begin{gathered}
\mathop {\sup }\limits_{\xi  \in [0,w]} | {f(\xi ,t) - \hat{h}(\xi
)} | \leqslant \tfrac{\sqrt{6}}{6}\theta w^{1/2} \text{ if } |f_\xi  (w,t)| > \theta , \hfill\\
\mathop {\sup }\limits_{\xi  \in [0,w]} |f(\xi ,t) - \hat{h}(\xi )
- 2 \theta w^{1/2}\sqrt {{\rm{v}}\,t} | \leqslant
\tfrac{\sqrt{6}}{6}\theta w^{1/2} \text{ if } |f_\xi (w,t)|
\leqslant \theta ,\hfill
\end{gathered}
\end{equation}
which completes the proof of Theorem~\ref{Th3}.
\end{proof}

{\footnotesize {\bf Acknowledgement.}  Author would like to thank
to Andreas  M\"{u}nch for his valuable comments and remarks.
Research is partially supported by the INTAS project Ref. No:
05-1000008-7921.}

\end{document}